\documentstyle[amssymb]{amsart}

\catcode`\@=11
\def\@setcopyright{}
\def\serieslogo@{}
\catcode`\@=12

\def\Q{{\bold Q}}
\def\Z{{\bold Z}}
\def\C{{\bold C}}
\def\R{{\bold R}}

\def\Gal{\mbox{Gal}}
\def\End{\mbox{End}}

\def\Hom{\mbox{Hom}}
\def\Lie{\mbox{Lie}}

\def\id{\mathrm{id}}
\def\GL{\mbox{GL}}
\def\Sp{\mbox{Sp}}
\def\U{\mbox{U}}
\def\SU{\mbox{SU}}
\def\dim{\mbox{dim}}

\def\det{\mbox{det}}
\def\tr{\mbox{tr}}
\def\O{{\cal O}}
\def\Pb{{\bold P}}

\def\N{{\bold N}}

\def\G{{\bold G}}
\def\Ss{{\bold S}}
\def\T{{\bold T}}

\def\Res{\mbox{Res}}
\def\SL{\mbox{SL}}
\def\zEnd{Z_{\mathrm{End}}}

\def\bmu{\hbox{\rm\rlap {$\mu$}\kern-.07em {$\mu$}}}

\newtheorem{thm}{Theorem}[section]
\newtheorem{lem}[thm]{Lemma}

\newtheorem{cor}[thm]{Corollary}
\newtheorem{prop}[thm]{Proposition}
\theoremstyle{definition}
\newtheorem{defn}[thm]{Definition}
\newtheorem{ex}[thm]{Example}
\newtheorem{rem}[thm]{Remark}
\title[Hodge groups of abelian varieties]
{Hodge groups of abelian varieties with purely multiplicative
reduction}

\author[A.\ Silverberg]{A.\ Silverberg}

\author[Yu.\ G.\ Zarhin]{Yu.\ G.\ Zarhin}

\begin{document}

\maketitle

\section{Introduction}

We show that if $A$ is an abelian variety over a subfield
$F$ of $\C$, and $A$ has purely multiplicative reduction at a discrete
valuation of $F$, then the Hodge group of $A$ is semisimple (Theorem
\ref{hodsemi}). Since the non-semisimplicity of the Hodge group of an
abelian variety can be translated into a condition on the endomorphism
algebra and its action on the tangent space (see Theorem \ref{notsemi}),
this gives a useful criterion for determining when an abelian variety
does not have purely multiplicative reduction.
For abelian varieties over number fields, a result analogous to
Theorem \ref{hodsemi} holds
where the Hodge group is replaced by a certain linear algebraic
group $H_\ell$ over $\Q_\ell$ arising from
the image of the $\ell$-adic representation associated to $A$ (see
Theorem \ref{notsemiGal}).
The Mumford-Tate conjecture predicts that $H_\ell$ is the
extension of scalars to $\Q_\ell$ of the Hodge group.

Our result generalizes a result of Mustafin (Corollary after
Theroem 3.2 of \cite{Mustafin}), which says that for a Hodge
family of abelian varieties (as in \cite{MumP}) admitting a
``strong degeneration'', generically the fibers have
semisimple Hodge group.
The problem of describing the Hodge group of an
abelian variety with
purely multiplication reduction was posed by V.\ G.\ Drinfeld,
in a conversation with Zarhin in the 1980's.

In \S\ref{boundssect} we provide bounds on torsion for
abelian varieties
which do not have purely multiplicative reduction at a
given discrete valuation. We apply this and Theorem \ref{hodsemi}
to obtain bounds on torsion for abelian varieties whose Hodge
groups are not
semisimple.

Silverberg would like to thank MSRI and IHES for their generous
hospitality, and NSF for financial support.
Zarhin would like to thank the Institute f\"ur Experimentelle
Mathematik for its hospitality,
and Gerhard Frey for his interest in the paper and useful
discussions.

\section{Definitions, notation, and lemmas}

Suppose $A$ is an abelian variety defined over a field $F$ of
characteristic zero, and $L$ is an algebraically
closed field containing $F$.
Write $\End_F(A)$
for the set of endomorphisms of $A$ which are defined over $F$, let
$\End(A) = \End_L(A)$, let $\End^0(A) = \End(A) \otimes_\Z \Q$,
and let $\End^0_F(A) = \End_F(A) \otimes_\Z \Q$.

Suppose $K$ is a field and $\iota : K \hookrightarrow \End_F^0(A)$
is an embedding such that $\iota(1) = 1$. Let $\Lie_F(A)$ be
the tangent
space of $A$ at the origin,
an $F$-vector space. If $\sigma$ is an embedding of $K$ into $L$, let
$$n_\sigma = \dim_L\{t \in \Lie_L(A) :
\iota(\alpha)t = \sigma(\alpha)t {\text{ for all }} \alpha \in K\}.$$
Note that $n_\sigma$ is independent of the choice of an algebraically
closed field $L$ containing $F$.
Write ${\bar \sigma}$ for the composition of
$\sigma$ with the involution complex conjugation of $K$.

\begin{defn}
If $A$ is an abelian variety over an algebraically
closed field $L$ of
characteristic zero, $K$ is a CM-field, and
$\iota : K \hookrightarrow \End^0(A)$
is an embedding such that $\iota(1) = 1$,
we say $(A,K,\iota)$ is {\em of Weil type}
if $n_\sigma = n_{\bar \sigma}$ for all embeddings $\sigma$
of $K$ into
$L$.
\end{defn}

\begin{lem}
\label{freeof}
If $A$ is an abelian variety defined over a field $F$ of characteristic
zero, $L$ is an algebraically closed field containing $F$,
$K$ is a CM-field,
and
$\iota : K \hookrightarrow \End_F^0(A) \subseteq \End^0(A)$
is an embedding such that $\iota(1) = 1$,
then the following statements are equivalent:
\begin{enumerate}
\item[\normalshape{(i)}]
$(A,K,\iota)$ is of Weil type,
\item[\normalshape{(ii)}]
$\iota$ makes $\Lie_{L}(A)$ into a free
$(K \otimes_\Q L)$-module,
\item[\normalshape{(iii)}]
$\iota$ makes $\Lie_F(A)$ into a free $(K \otimes_\Q F)$-module.
\end{enumerate}
\end{lem}

\begin{pf}
Let
$\Sigma$ be the set of embeddings of $K$ into $L$, let
$\psi_F : K \otimes_\Q F \to \End(\Lie_F(A))$
be the homomorphism induced by $\iota$,
let $\psi : K \otimes_\Q L \to \End(\Lie_L(A))$ be the extension
of scalars of $\psi_F$ to $K \otimes_\Q L$, let
$m = 2\dim(A)/[K : \Q]$, let $M_F = (K \otimes_\Q F)^m$, and let
$M = (K \otimes_\Q L)^m = M_F \otimes_F L$. Let
$\psi_F^\prime : K \otimes_\Q F \to \End(M_F)$
and
$\psi^\prime : K \otimes_\Q L \to \End(M)$ be the natural homomorphisms.
By \S 2.1 of \cite{Shimura}, for every $\sigma \in \Sigma$ we have
$n_\sigma + n_{\bar \sigma} = m$.
For $\alpha \in K$, taking the trace of $\psi(\alpha)$ gives
$$\tr(\psi(\alpha)) = \sum_{\sigma \in \Sigma}n_{\sigma}\sigma(\alpha).$$
The traces of $\psi$ and of $\psi^\prime$ coincide on $K$ if and only if
$n_\sigma = n_{\bar \sigma} = m/2$ for every $\sigma \in \Sigma$.
Since $K \otimes_\Q L$ is a semisimple ring,
$\Lie_L(A)$ and $M$ are semisimple $(K \otimes_\Q L)$-modules.
Therefore, $\Lie_L(A)$ is a free $(K \otimes_\Q L)$-module if and only if
the traces of $\psi$ and of $\psi^\prime$ coincide on $K$.
Therefore, $\Lie_L(A)$ is a free $(K \otimes_\Q L)$-module if and only if
$(A,K,\iota)$ is of Weil type.

If $\Lie_F(A)$ is a free
$(K \otimes_\Q F)$-module, then clearly
$\Lie_L(A)$ ( $= \Lie_F(A) \otimes_F L$) is a free
$(K \otimes_\Q L)$-module.

Conversely, if $\Lie_L(A)$ is a free
$(K \otimes_\Q L)$-module, then the traces of $\psi$ and of
$\psi^\prime$ (and therefore of $\psi_F$ and of
$\psi_F^\prime$) coincide on $K$. Since
$K \otimes_\Q F$ is a semisimple ring, $\Lie_F(A)$ and $M_F$ are
semisimple $(K \otimes_\Q F)$-modules. Therefore $\Lie_F(A)$ and $M_F$
are isomorphic as $(K \otimes_\Q F)$-modules, i.e., $\Lie_F(A)$
is a free $(K \otimes_\Q F)$-module.
\end{pf}

See also p.~525 of \cite{Ribet} for
the case where $K$ is an imaginary quadratic field.

\begin{rem}
If $A$ is an abelian variety defined over a field $F$ of characteristic
zero, $K$ is a totally real number field, and $\iota : K \hookrightarrow
\End_F^0(A)$ is an embedding such that $\iota(1) = 1$, then
$\iota$ makes $\Lie_F(A)$ into a free $(K \otimes_\Q F)$-module.
To see this, let $L$ be an algebraically closed field containing $F$, let
$\Sigma$ be the set of embeddings of $K$ into $L$, and let
$\psi : K \to \End(\Lie_L(A))$
be the homomorphism induced by $\iota$.
Let $m^{\prime} = \dim(A)/[K : \Q]$. We have
$$\tr(\psi(\alpha)) = m^{\prime}\sum_{\sigma \in \Sigma}\sigma(\alpha)$$
for every $\alpha \in K$, by \S 2.1 of \cite{Shimura}. Therefore,
$\Lie_L(A)$ is a free $(K \otimes_\Q L)$-module. As in the
proof of Lemma \ref{freeof}, it follows that $\Lie_F(A)$ is a free
$(K \otimes_\Q F)$-module.
\end{rem}

Suppose $A$ is a complex abelian variety.
Let $V = H_1(A,\Q)$ and let $\Ss = \Res_{\C/\R}\G_m$.
The complex structure on $A$ gives rise to a rational Hodge structure
on $V$ of weight $-1$, i.e., a homomorphism of algebraic groups
$h : \Ss \to \GL(V)_\R$.
Let $\T$ be the kernel of the norm map $\N : \Ss \to \G_m$.
Then $\T(\R) = \{x \in \C : |x| = 1\}$.

\begin{defn}
If $A$ is an abelian variety over $\C$ and $V = H_1(A,\Q)$, then the
{\em Hodge group} $H$ is the smallest
algebraic subgroup of $\GL(V)$ defined over $\Q$ such that $H(\R)$
contains $h(\T(\R))$.
Equivalently, $H$ is the largest algebraic subgroup of $\GL(V)$
defined over $\Q$ such that all Hodge classes in
$V^{\otimes p} \otimes (V^\ast)^{\otimes q}$, for all non-negative
integers $p$ and $q$, are tensor invariants of $H$. I.e.,
$H$ is the largest algebraic subgroup of $\GL(V)$ defined over $\Q$
which fixes all Hodge classes of all powers of $A$.
\end{defn}

It follows from the definition of $H$ that $\End^0(A) = \End_H(V)$.

If now $F$ is a number field and $\ell$ is a prime number, let
$T_\ell(A) = {\displaystyle \lim_\leftarrow A_{\ell^r}}$
(the Tate module), let
$V_\ell(A) = T_\ell(A) \otimes_{\Z_\ell}\Q_\ell$, and let $\rho_{A,\ell}$
denote the $\ell$-adic representation
$$\rho_{A,\ell} : \Gal({\bar F}/F) \to \GL(T_\ell(A))
\subseteq \GL(V_\ell(A)).$$
Let $G_\ell$
denote the algebraic envelope of the image of $\rho_{A,\ell}$, i.e., the
Zariski closure in $\GL(V_\ell(A))$ of the image of $\rho_{A,\ell}$.
By \cite{Faltings}, $G_\ell$ is a reductive algebraic group, and
$\End^0_F(A) \otimes_\Q \Q_\ell  = \End_{G_\ell}(V_\ell(A))$.
Let $H_\ell$ be the identity connected component of
$G_\ell \cap \SL(V_\ell(A))$.
Then $H_\ell$ is
a connected reductive group and
$\End^0(A) \otimes_\Q \Q_\ell = \End_{H_\ell}(V_\ell(A))$.

We will repeatedly use the fact (see the first Theorem on p.~220
of \cite{Hum})
that if $G$ is a connected linear algebraic group over a field $F$
of characteristic zero, then $G(F)$ is Zariski-dense in $G$.

\begin{lem}
\label{semfin}
If $G$ is a reductive linear algebraic group over a
field $F$ of characteristic zero, and $Z$ is the center of $G$,
then $G$ is semisimple if and only if $Z(F)$ is finite.
\end{lem}

\begin{pf}
Let $Z^0$ denote the identity connected component of $Z$.
Since $G$ is reductive, $G$ is semisimple if and only if $Z^0 = 1$
(see the lemma on p.~125 of \cite{Hum}). Since $Z^0(F)$ is Zariski-dense
in $Z^0$, $Z^0 = 1$ if and only if $Z(F)$ is finite.
\end{pf}

\section{Semisimplicity criteria for the groups $H$ and $H_\ell$}

If the center of $\End^0(A)$
is a direct sum of totally real number fields, then
it is well-known that
the groups $H$ and $H_\ell$ are semisimple (see, for instance,
Corollary 1 in \S 1.3.1 of
\cite{Zarhintor} and Lemma 1.4 of \cite{Tankeev}).
The following result follows easily from a result in \cite{moonzar}, and
characterizes the endomorphism algebras of abelian varieties whose
Hodge groups are not semisimple.

\begin{thm}
\label{notsemi}
Suppose $A$ is an
abelian variety defined over $\C$.
Then the Hodge group of $A$ is not semisimple if and only if
for some simple component $B$ of $A$, the
center of $\End^0(B)$ is a CM-field $K$ such that
$(B,K,\id)$ is not of Weil type, with $\id$ the identity embedding
of $K$ in $\End^0(B)$.
\end{thm}

\begin{pf}
Let $V = H_1(A,\Q)$.
Fix a polarization on $A$. The polarization induces a
non-degenerate alternating bilinear form
$\varphi : V \times V \to \Q$ such that
$H \subseteq \Sp(V,\varphi)$ (see \cite{MumMatAnn}). Then
$$H(\Q) \subseteq \Sp(V,\varphi)(\Q) = \{ g \in \End(V) : gg^\prime = 1 \},$$
where $g \mapsto g^\prime$ is the involution on $\End(V)$ defined
by
$$\varphi(g(x),y) = \varphi(x,g^\prime (y)) \mbox{ for $x, y \in V$}.$$
The restriction of the involution ${}^\prime$ to $\End^0(A)$
is the Rosati involution.
Let $Z$ denote the center of $H$ and let $\zEnd$ denote the center of
$\End^0(A)$.
If $\alpha \in Z(\Q)$, then $\alpha$ commutes with all elements of
$H(\Q)$, so $\alpha \in \End^0(A)$. Further, since $\alpha \in H(\Q)$,
$\alpha$ commutes with all elements of $\End^0(A)$, and therefore
$\alpha \in \zEnd$. Therefore,
\begin{equation}
\label{Zs}
Z(\Q) \subseteq \{\alpha \in \zEnd : \alpha\alpha^{\prime} = 1\}.
\end{equation}
If $A$ is isogenous to a product of two abelian varieties, then
the Hodge group $H$ of $A$ is a
subgroup of the product of the Hodge groups $H_1$ and $H_2$ of
the factors, in such a way that for $i = 1$ and $2$ the restriction to $H$
of the projection map from $H_1 \times H_2$ onto $H_i$ induces a
surjective homomorphism from $H$ onto $H_i$ (see Proposition 1.6 of
\cite{Hazama}). It follows easily that
$H$ is semisimple if and only if both $H_1$ and $H_2$ are semisimple.
We may therefore reduce to the case where $A$ is a simple abelian
variety.
Then the center $\zEnd$ of $\End^0(A)$ is either a totally real number
field or a CM-field.

Suppose $\zEnd$ is totally real. Then all Rosati
involutions are the identity when restricted to $\zEnd$.
By (\ref{Zs}), $Z(\Q) \subseteq \{\pm 1\}$. Therefore,
$Z(\Q)$ is finite, so $H$ is semisimple by Lemma \ref{semfin}.

Suppose $\zEnd$ is a CM-field $K$.
Then every Rosati involution induces complex conjugation on $K$.
Choose $\alpha \in K^\times$ such that ${\bar \alpha} = -\alpha$.
Then there exists a unique $K$-Hermitian form $\psi : V \times V \to K$
such that $\varphi(x,y) = Tr_{K/\Q}(\alpha \psi(x,y))$ (see \cite{Shimura}).
The unitary group $\U(V,\psi)$ is an algebraic group over $K_0$, the
maximal totally real subfield of $K$. Let
$U = \Res_{K_0/\Q} \U(V,\psi)$,
let $\SU$ denote the kernel of the determinant homomorphism
$\det_K : \U \to \Res_{K/\Q} \G_m$, and let
$\End_K(V)$ denote the ring of $K$-linear endomorphisms of $V$.
Then
$$\U(\Q) = \{g \in \End_K(V) :
\psi(g(x),g(y)) = \psi(x,y)  {\text{ for all }}  x, y \in V\}$$
and $H \subseteq \U \subseteq \Sp(V,\varphi)$.
By Lemma 2.8 of \cite{moonzar}, $H \subseteq \SU$
if and only if $(A,K,\id)$ is of Weil type.
If $H$ is semisimple, then all homomorphisms from $H$ to
commutative groups are trivial. Therefore
$\det_K(H) = 1$, so $H \subseteq \SU$.
Conversely,
if $H \subseteq \SU$, then
$Z(\Q) \subseteq \SU(\Q) \cap K$,
the group of $(\dim(V)/[K:\Q])$-th
roots of unity in $K$. Therefore $Z(\Q)$ is
finite and $H$ is semisimple.
\end{pf}

\begin{thm}
\label{notsemiGal}
Suppose $A$ is an abelian variety defined over a number field $F$. Then
the following are equivalent:
\begin{enumerate}
\item[\normalshape{(i)}] $H$ is semisimple,
\item[\normalshape{(ii)}] $H_\ell$ is semisimple, for one prime $\ell$,
\item[\normalshape{(iii)}] $H_\ell$ is semisimple, for every prime $\ell$.
\end{enumerate}
\end{thm}

\begin{pf}
Let $\ell$ be a prime number and let $V_\ell = V_\ell(A)$.
By Theorem \ref{notsemi}, it suffices to show that
$H_\ell$ is not
semisimple if and only if for some simple component $B$ of $A$, the
center of $\End^0(B)$ is a CM-field $K$ such that
$(B,K,\id)$ is not of Weil type, with $\id$ the identity embedding
of $K$ in $\End^0(B)$.
Since $H_\ell$ is connected, it is invariant under
finite extensions of the number field $F$.
By replacing $F$ by a finite extension, we may suppose
that $\End^0(A) = \End_F^0(A)$.

We parallel the proof of Theorem \ref{notsemi}. Fix a polarization
on $A$ defined over $F$.
Let $V$ and $\varphi$ be as in the proof of Theorem \ref{notsemi}.
Then $V_\ell = V \otimes_\Q \Q_\ell$.
Let
$\varphi_\ell : V_\ell \times V_\ell \to \Q_\ell$ be the $\Q_\ell$-linear
extension of $\varphi$.
It follows
immediately from p.~516 of \cite{zarhin} and the definition of $H_\ell$ that
$H_\ell \subseteq \Sp(V_\ell,\varphi_\ell)$.
Let $Z_\ell$ denote the center of $H_\ell$, let $\zEnd$
denote the center of $\End^0(A)$, and let
${}^\prime$ denote the involution on $\End(V_\ell)$ induced by
$\varphi_\ell$. Following the proof of
Theorem \ref{notsemi}, we conclude that
$$Z_\ell(\Q_\ell) \subseteq \{\alpha \in \zEnd \otimes_\Q \Q_\ell :
\alpha\alpha^{\prime} = 1\}.$$
If $A$ is $F$-isogenous to a product of two abelian varieties, then
$H_\ell$ is a
subgroup of the product of the corresponding groups $H_{1,\ell}$ and
$H_{2,\ell}$ for
the factors, in such a way that for $i = 1$ and $2$ the restriction to $H_\ell$
of the projection map from $H_{1,\ell} \times H_{2,\ell}$ onto
$H_{i,\ell}$ induces a
surjective homomorphism from $H_\ell$ onto $H_{i,\ell}$. It follows that
we may reduce to the case where $A$ is $F$-simple.

If $\zEnd$ is totally real, we conclude that $H_\ell$ is semisimple
as in the proof of Theorem \ref{notsemi}.
Suppose $\zEnd$ is a CM-field $K$ and
let $K_\ell = K \otimes_\Q \Q_\ell$. Let $\psi$ and $U$ be as in
the proof of Theorem \ref{notsemi},
let  $\psi_\ell : V_\ell \times V_\ell \to K_\ell$
denote the $K_\ell$-Hermitian form which extends the pairing $\psi$,
let $\U_\ell = \U \times \Q_\ell$,
let $\SU_\ell$ denote the kernel of the determinant homomorphism
$\det_{K_\ell} : \U_\ell \to \Res_{K/\Q}\G_m \times \Q_\ell$, and let
$\End_{K_\ell}(V_\ell)$ denote the ring of
$K_\ell$-linear endomorphisms of $V_\ell$. Then
$$\U_\ell(\Q_\ell) =
\{g \in \End_{K_\ell}(V_\ell) :
\psi_\ell(g(x),g(y)) = \psi_\ell(x,y)  {\text{ for all }}  x, y \in V_\ell\}$$
and
$H_\ell \subseteq \U_\ell \subseteq \Sp(V_\ell,\varphi_\ell)$.
By Lemma 2.8 of \cite{moonzar}, $H_\ell \subseteq \SU_\ell$
if and only if $(A,K,\id)$ is of Weil type.
The group  $\SU_\ell(\Q_\ell) \cap K_\ell$ is
the finite group of
$(\dim_{\Q_\ell}(V_\ell)/\dim_{\Q_\ell}(K_\ell))$-th
roots of unity in the ring $K_\ell$.
Paralleling the proof of Theorem \ref{notsemi},
$H_\ell \subseteq \SU_\ell$
if and only if $H_\ell$ is semisimple.
\end{pf}

\begin{ex}
If $A$ is odd-dimensional and the center of
$\End^0(A)$ is a CM-field $K$, then $H$ is not semisimple.
To show this, note that $A$ is isogenous to a power of a simple
odd-dimensional abelian variety $B$ such that $K$ is the center
of $\End^0(B)$. Then the Hodge groups of $A$ and of $B$ coincide,
so we may reduce to the case where $A$ is simple. Let
$d = \dim(A)$ and use the
notation of the proof of Lemma \ref{freeof}. Then
$n_\sigma +  n_{\bar \sigma} = 2d/[K : \Q]$.
If $H$ were semisimple, then by Theorem \ref{notsemi},
$(A,K,\id)$ would be of Weil type. We would therefore have
$n_\sigma = n_{\bar \sigma}$, and so $2d/[K : \Q]$ would be even.
However, $d$ is odd and $[K : \Q]$ is even, so this cannot happen.
\end{ex}

\section{Abelian varieties having purely multiplicative reduction}

\begin{thm}
\label{hodsemi}
Suppose $A$ is an abelian variety over a subfield $F$ of $\C$,
$v$ is a discrete valuation on $F$, and $A$ has purely
multiplicative reduction
at $v$. Then the Hodge group $H$ of $A$ is semisimple.
\end{thm}

\begin{pf}
Since $H$ is semisimple if the Hodge groups of each of its
$F$-simple components are,
we may reduce to the case where $A$ is an $F$-simple abelian
variety with purely multiplicative reduction at $v$.
Since the properties of having semisimple Hodge group and having
purely multiplicative reduction are invariant under finite
extensions of the ground field, we may assume
$\End_F^0(A) = \End^0(A)$.
Suppose that $H$ is not semisimple.
By Theorem \ref{notsemi}, the
center of $\End^0(A)$ is a CM-field $K$ such that $(A,K,\id)$ is not
of Weil type, with $\id$ the identity embedding
of $K$ in $\End^0(A)$.
Let $L$ be a fixed algebraic closure of the completion $F_v$ of
$F$ at $v$. Since $A$ has purely multiplicative reduction at $v$,
$A$ admits a non-archimedean uniformation; i.e.,
(see \cite{Mumford} and \cite{Raynaud}) there are a discrete
subgroup $\Gamma$ of $\G_m^d(F_v) = (F_v^{\times})^d$,
isomorphic to $\Z^d$, and a $\Gal(L/F_v)$-equivariant $v$-adically
continuous isomorphism $(L^{\times})^d/\Gamma \cong A(L)$ which for
some finite extension $M$ of $F_v$ induces
an isomorphism $(M^{\times})^d/\Gamma \cong A(M)$
as $M$-Lie groups.
Let $\O$ be the center of $\End(A)$. Then $\O$ is an order in $K$.
By Satz 6 of \cite{Gerritzen}, there is a homomorphism $\O \hookrightarrow
\End(\G_m^d)$ which induces the inclusion $\O \subseteq \End(A)$.
Composing with the natural homomorphism
$$\End(\G_m^d) \hookrightarrow \End(\Hom(\G_m,\G_m^d))$$
and tensoring with $\Q$, we have
$$K \hookrightarrow \End(\Hom(\G_m,\G_m^d) \otimes \Q).$$
Therefore, the inclusion of $K$ in $\End^0(A)$ induces a $K$-vector space
structure on $\Hom(\G_m,\G_m^d) \otimes \Q$.
Tensoring with $M$ makes $\Hom(\G_m,\G_m^d) \otimes_\Z M$ ( $ = M^d$) into
a free $(K \otimes_\Q M)$-module. We can view $(M^\times)^d$
as a (non-archimedean) analytic variety over $M$. The tangent
space to
$(M^{\times})^d$ at $1$ is isomorphic to $M^d$.
By \cite{Morik1} (see also Chapter 2 of \cite{Manin}),
$(L^{\times})^d/\Gamma$ can be embedded, via
theta functions, as an analytic subvariety of a projective space
$\Pb^n(L)$, so that the image of $(M^{\times})^d/\Gamma$ is $A(M)$.
Let $T$ denote the analytic tangent space at the origin of the
analytic variety $A(M)$. The tangent map is an isomorphism
$M^d \cong T$.
The algebraic tangent space at the origin to the algebraic variety
$A$ over $M$ is $\Lie_M(A) = \Lie_F(A) \otimes_F M$, and there is
a canonical isomorphism between the analytic and algebraic
tangent spaces to $A(M)$
(see subsection 3 of \S 2 of Chapter II of \cite{shaf}). Therefore,
the identity embedding of $K$ into $\End^0(A)$ makes
$\Lie_M(A)$ into a free $(K \otimes_\Q M)$-module.
By Lemma \ref{freeof}, $(A,K,\id)$
is of Weil type, contradicting our assumptions.
\end{pf}

Theorem \ref{hodsemi} remains true if we replace the assumption
that $v$ is a discrete valuation by the assumption that $v$ is a
valuation of rank $1$ and $A$ admits non-archimedean
uniformization (Gerritzen's theorem remains true under these assumptions).

\section
{Bounds on torsion of abelian varieties which do not have purely
multiplicative reduction}
\label{boundssect}

It is easy to find uniform bounds on orders of torsion points over
number fields for
abelian varieties with potential good reduction, or for elliptic
curves which do not have multiplicative reduction, at a given
discrete valuation (see \cite{contemp}, \cite{CM}, \cite{flexoest}).
In this section we extend these results by finding bounds on torsion
subgroups of abelian varieties
which do not have purely multiplicative reduction at a given
discrete valuation.

Suppose $A$ is a $d$-dimensional abelian variety over a field $F$,
$v$ is a discrete valuation on $F$ with finite residue field $k$
of order $q$, $n$ is a positive integer
relatively prime to $q$, and $J$ is a non-zero subgroup of the group
$A_n(F)$ of points in $A(F)$ of order dividing $n$.
Let $A_v^0$ denote the connected component of the identity of the
special fiber $A_v$ of the N\'eron minimal model of $A$ at $v$.
Let $a$, $u$, and $t$ denote respectively the abelian, unipotent,
and toric ranks of $A_v^0$. Then $d = a + t + u$.

If $\lambda$ is a polarization on $A$
defined over an extension of $F$ which is unramified over $v$,
define a skew-symmetric Galois-equivariant pairing
$e_{\lambda,n}$ on $A_n$ by
$e_{\lambda,n}(x,y) = e_n(x,\lambda(y))$, where $e_n$ is the
Weil pairing.
If $J$ is not isotropic
with respect to $e_{\lambda,n}$, and $n$ is a prime number,
then $\zeta_n \in F$, so the prime $n$ can be bounded independent
of $A$ (with a bound depending on $F$).
Therefore, the more interesting case is when $J$ is an isotropic
subgroup of $A_n(F)$.
If $J$ is a maximal isotropic subgroup of $A_n(F)$, and
$A$ does not have semistable reduction at $v$, then $n \le 4$,
by Theorem 6.2 of \cite{semistab}.
The remaining case to consider
is the case where $A$ has semistable reduction at $v$. Theorem \ref{bounds}
below implies that in
this case we can bound $n$ in terms of $q$ and $d$, as long
as $A$ does not have purely multiplicative reduction at $v$.
Note that if $P$ is a point of $A(F)$ of order $n$ which reduces to a point of
$A_v^0$, then $n$ is bounded above by $\#A_v^0(k)$. Therefore
even in the case of abelian varieties with purely
multiplicative reduction one can easily bound, by
a constant depending only on $d$ and $q$,
the orders of torsion points
whose reductions lie in $A_v^0$. As was the case for
elliptic curves, the most difficult case is the case when the
reduction is purely multiplicative and the reductions of the torsion
points do not lie in the identity connected component of the
special fiber of the N\'eron minimal model.

\begin{lem}[Lemma 1 on pp.~494--495 of \cite{serretate}]
\label{serretatelem}
If $A$ is an abelian variety over a field $F$,
$v$ is a discrete valuation on $F$, and $n$ is a positive
integer relatively prime to the residue characteristic of $v$,
then $(A_v^0)_n$ is a free $\Z/n\Z$-module of rank $2a + t$.
\end{lem}

\begin{prop}
\label{lessthan}
Suppose $A$ is an abelian variety over a field $F$,
$v$ is a discrete valuation on $F$ with finite residue field $k$
of order $q$, $n$ is a positive integer relatively prime to $q$,
and $J$ is a subgroup of $A_n(F)$.
Suppose there is a positive constant $\epsilon$ such that
$|J| \ge n^{t+2u+\epsilon}$.
Then
$n \le (\#(A_v^0)_n(k))^{1/\epsilon} \le (\#A_v^0(k))^{1/\epsilon}$.
\end{prop}

\begin{pf}
Let $d = \dim(A)$.
Via the reduction map we may view $A_n(F)$, and therefore $J$,
as a subgroup of $(A_v)_n$ (see \cite{serretate}).
Therefore,  $\#J \#(A_v^0)_n$ divides $n^{2d} \#(J \cap (A_v^0)_n)$.
Thus by Lemma \ref{serretatelem},
$\#J$ divides $n^{t+2u}\#(J \cap (A_v^0)_n)$.
Therefore,
$$n^\epsilon \le \#(J \cap (A_v^0)_n) \le
\#(A_v^0)_n(k) \le \#A_v^0(k).$$
\end{pf}

\begin{thm}
\label{bounds}
Suppose $A$ is a $d$-dimensional abelian variety over a field $F$,
$v$ is a discrete valuation on $F$ with finite residue field $k$
of order $q$, $n$ is a positive integer relatively prime to
$q$, and $A_n(F)$ has a subgroup of
order $n^d$. Suppose the reduction of $A$ at $v$ is semistable but
not purely multiplicative.
Then
$$n \le (1 + \sqrt{q})^{2a}(1 + q)^t \le (1 + \sqrt{q})^{2d}.$$
\end{thm}

\begin{pf}
Since $A$ has semistable reduction at $v$, $u = 0$.
Since the reduction of $A$ at $v$
is not purely multiplicative, $a \ge 1$.
Applying Proposition \ref{lessthan} with
$\epsilon = 1$, we have $n \le \#A_v^0(k)$.
Since $A$ has semistable reduction at $v$, $A_v^0$ is an extension
of an abelian variety $B$ by a torus $T$. We have the Weil bound
$\# B(k) \le (1 + \sqrt{q})^{2a}$. Similarly, we have the bound
$\# T(k) \le (1 + q)^t$, as follows.
Let $X$ be the group of characters of $T \otimes {\bar k}$.
The Frobenius element of $\Gal({\bar k}/k)$ acts on $X$, say by
$\varphi_0$. Since the torus $T$
splits over some finite extension of $k$, $\Gal({\bar k}/k)$ acts on
$X$ through a finite quotient, so
all the eigenvalues of $\varphi_0$ have absolute value $1$.
Therefore all eigenvalues of $q - \varphi_0$ are non-zero and
have absolute value at most $1 + q$. We have (see Theorem 6.2
in \S 1 of Chapter VI of \cite{Vo})
$$\# T(k) = |\det(q - \varphi_0)| \le (1 + q)^t.$$
Therefore,
$$\#A_v^0(k) \le  (1 + \sqrt{q})^{2a}(1 + q)^t =
(1 + \sqrt{q})^{2a}(1 + q)^{d-a}
\le (1 + \sqrt{q})^{2d}.$$
\end{pf}

If $c$ and $d$ are positive integers, let $f(c,d)$ be the maximum of
the orders of the elements of $GL_{2d}(\Z/c\Z)$.

\begin{thm}
\label{nrbounds}
Suppose $A$ is a $d$-dimensional abelian variety over a number field $F$
of degree $m$, $v$ is a discrete valuation on $F$ at which $A$ does not
have purely multiplicative reduction, $p$ is the residue
characteristic of $v$, $n$ and $r$ are positive integers not divisible
by $p$, $r \ge 3$, and $A_n(F)$ has a subgroup of order $n^d$. Then
$$n \le (1 + p^{mf(r,d)/2})^{2d}.$$
\end{thm}

\begin{pf}
By a theorem of Raynaud (Proposition 4.7 of \cite{SGA}), $A$ has semistable
reduction at the discrete valuations on the field $F(A_{r})$ of residue
characteristic not dividing $r$. Let $v^\prime$ be a valuation
on $F(A_r)$ extending $v$, and let $k$ be the corresponding residue
field. Then $\# k$ divides $p^{mf}$, where $f$ is the order of
Frobenius at $v^\prime$ in $\Gal(F(A_r)/F)$. Since
$\Gal(F(A_r)/F)$ injects into $GL_{2d}(\Z/r\Z)$,
$\# k$ divides $p^{mf(r,d)}$. The result now follows from Theorem
\ref{bounds}.
\end{pf}

\begin{cor}
\label{Hsemibds}
Suppose $A$ is a $d$-dimensional abelian variety over a number field $F$
of degree $m$, and suppose the Hodge group of $A$ is not semisimple.
Suppose $n$ is a positive integer and $A_n(F)$ has a subgroup of
order $n^d$. Then
$$n \le [(1 + 2^{mf(3,d)/2})(1 + 3^{mf(4,d)/2})]^{2d}
< (1+10^{-11})\cdot (2^{3^{4d^2}}\cdot 3^{4^{4d^2}})^{md}.$$
\end{cor}

\begin{pf}
The result follows from Theorem \ref{hodsemi}, by applying
Theorem \ref{nrbounds} with
$p = 2$, $r = 3$ to bound the prime-to-two part of $n$,
and with $p = 3$, $r = 4$ to bound the prime-to-three part of $n$.
The final inequality follows from the bound
$f(c,d) \le \# GL_{2d}(\Z/c\Z) < c^{4d^2}$.
\end{pf}

The bounds on $n$ given in Corollary \ref{Hsemibds} were shown
in Theorem 3.3 and Remark 2 of \cite{contemp} to be bounds
on the orders of torsion subgroups of
abelian varieties with potential good reduction at discrete valuations
of residue characteristics $2$ and $3$.

If we assume the existence of a polarization on $A$ of degree prime
to $n$ (for example, a principal polarization) we obtain stronger
bounds. The following results give bounds on torsion subgroups of
order prime to the degree of a given polarization.

\begin{thm}
\label{lambdabounds}
Suppose $A$ is a $d$-dimensional abelian variety over a number field $F$
of degree $m$, $v$ is a discrete valuation on $F$ at which $A$ does not
have purely multiplicative reduction, $p$ is the residue
characteristic of $v$, $\ell$ is a prime number, $\ell \ne p$, $J$ is
a subgroup of $A_\ell(F)$ of order $\ell^d$, $\lambda$ is a polarization
on $A$ defined over an extension of $F$ unramified at $v$, and $\ell$
does not divide the degree of $\lambda$. Then
$$\ell \le (1 + p^{m/2})^{2d}.$$
\end{thm}

\begin{pf}
Since $\ell$ does not divide the degree of $\lambda$,
the pairing $e_{\lambda,\ell}$ is nondegenerate.
If $J$ is not isotropic with respect to $e_{\lambda,\ell}$, then $F$
contains a primitive $\ell$-th root of unity. Therefore $\ell - 1$
divides $[F : \Q]$, so $\ell \le 1 + m$.

Suppose $J$ is isotropic with respect to $e_{\lambda,\ell}$.
Then $J$ is a maximal isotropic subgroup of $A_\ell$ (since
$\# J = \ell^d$ and $e_{\lambda,\ell}$ is nondegenerate).
If $A$ does not have semistable reduction at $v$, then $\ell \le 3$
by Theorem 6.2 of \cite{semistab}. If $A$ has semistable reduction
at $v$, then $\ell \le (1 + p^{m/2})^{2d}$ by Theorem \ref{bounds}.
The result now follows since $(1 + p^{m/2})^{2d}$ is greater than
$3$ and than $1 + m$.
\end{pf}

\begin{cor}
\label{Hsemilambdabds}
Suppose $(A,\lambda)$ is a $d$-dimensional polarized abelian variety
over a number field $F$ of degree $m$, and suppose the Hodge group of
$A$ is not semisimple.
Suppose $\ell$ is a prime number which does not divide the degree of
$\lambda$, and $J$ is a subgroup of $A_\ell(F)$ of order $\ell^d$.
Then
$$\ell \le (1 + 2^{m/2})^{2d}.$$
\end{cor}

\begin{pf}
By  Theorem \ref{hodsemi}, $A$ does not have
purely multiplicative reduction at any discrete valuations.
Since $2 < (1 + 2^{m/2})^{2d}$, we may assume $\ell$ is an odd prime, and
we obtain the result by applying Theorem \ref{lambdabounds} with $p = 2$.
\end{pf}

The proof of Theorem \ref{lambdabounds} shows that
$\ell \le \max\{ 1+m, (1+\sqrt{q})^{2d}\}$, where $q$ is the order
of the residue field of $v$. Therefore in Corollary \ref{Hsemilambdabds}
we can conclude that $\ell \le \max\{ 1+m, (1+\sqrt{f})^{2d}\}$, where $f$
is the minimal order of the residue fields of the valuations on $F$ of
residue characteristic $2$.

\begin{cor}
Suppose $(A,\lambda)$ is a $d$-dimensional polarized abelian variety
over a number field $F$ of degree $m$, and suppose the Hodge group of
$A$ is not semisimple.
Suppose $n$ is a positive integer relatively prime to the degree of
$\lambda$, and $J$ is a subgroup of $A_n(F)$ of order $n^d$
which is a maximal isotropic subgroup with
respect to $e_{\lambda,n}$.
Then
$$n \le (1 + 2^{m/2})^{2d}(1 + 3^{m/2})^{2d}.$$
\end{cor}

\begin{pf}
By  Theorem \ref{hodsemi}, $A$ does not have
purely multiplicative reduction at any discrete valuations.
The prime-to-$p$ part of $n$ is
bounded above by $4$ if there is a valuation on $F$ of
residue characteristic
$p$ at which $A$ does not have semistable reduction (by Theorem 6.2 of
\cite{semistab}), and otherwise is bounded above by $(1 + p^{m/2})^{2d}$
(by Theorem \ref{bounds}). Note that $4 < (1 + p^{m/2})^{2d}$.
The result follows by letting $p = 2, 3$.
\end{pf}

\bigskip

\noindent {\small {Mathematics Department,}}

\noindent {\small {Ohio State University,}}

\noindent {\small {Columbus, OH 43210-1174,}}

\noindent {\small {USA}}

\noindent {\small {\em E-mail address}: silver@math.ohio-state.edu}

\bigskip

\noindent {\small {Mathematics Department,}}

\noindent {\small {Pennsylvania State University,}}

\noindent {\small {University Park, PA 16802,}}

\noindent {\small {USA} }

\noindent and

\noindent {\small {Institute for Mathematical Problems in Biology, }}

\noindent {\small {Russian Academy of Sciences, }}

\noindent {\small {Push\-chino, Moscow Region, 142292, }}

\noindent {\small {RUSSIA}}

\noindent {\small {\em E-mail address}: zarhin@math.psu.edu}

\end{document}